
\def\LaTeX{%
  \let\Begin\begin
  \let\End\end
  \let\salta\relax
  \let\finqui\relax
  \let\futuro\relax}

\def\UK{\def\our{our}\let\sz s}
\def\USA{\def\our{or}\let\sz z}



\LaTeX

\USA


\salta

\documentclass[twoside,12pt]{article}
\setlength{\textheight}{24cm}
\setlength{\textwidth}{16cm}
\setlength{\oddsidemargin}{2mm}
\setlength{\evensidemargin}{2mm}
\setlength{\topmargin}{-15mm}
\parskip2mm


\usepackage{amsmath}
\usepackage{amsthm}
\usepackage{amssymb}
\usepackage[mathcal]{euscript}

\usepackage[usenames,dvipsnames]{color}
%
%


\def\gianni{\color{red}}
\def\Gianni{\color{green}}
\def\GP{\color{blue}}
\let\gianni\relax
\let\Gianni\relax
\let\GP\relax

\def\blu #1{\textcolor{blue}{#1}}
\definecolor{viol}{rgb}{0.5,0,0.3}
\def\vio #1{\textcolor{viol}{#1}}
\let\blu\relax
\let\vio\relax
%

\bibliographystyle{plain}


%

\finqui

\def\Beq{\Begin{equation}}
\def\Eeq{\End{equation}}
\def\Bsist{\Begin{eqnarray}}
\def\Esist{\End{eqnarray}}

\def\Bthm{\Begin{theorem}}
\def\Ethm{\End{theorem}}
\def\Blem{\Begin{lemma}}
\def\Elem{\End{lemma}}
\def\Bprop{\Begin{proposition}}
\def\Eprop{\End{proposition}}
\def\Bcor{\Begin{corollary}}
\def\Ecor{\End{corollary}}
\def\Brem{\Begin{remark}\rm}
\def\Erem{\End{remark}}

\def\Bdim{\Begin{proof}}
\def\Edim{\End{proof}}
\let\non\nonumber




\def\step #1 \par{\medskip\noindent{\bf #1.}\quad}


\def\Lip{Lip\-schitz}
\def\holder{H\"older}
\def\aand{\quad\hbox{and}\quad}

\def\lhs{left-hand side}
\def\rhs{right-hand side}
\def\sfw{straightforward}


\def\generaliz{generali\sz}

\def\organiz{organi\sz}

\def\summariz{summari\sz}

\def\bhv{behavi\our}


\def\multibold #1{\def\arg{#1}%
  \ifx\arg\pto \let\next\relax
  \else
  \def\next{\expandafter
    \def\csname #1#1#1\endcsname{{\bf #1}}%
    \multibold}%
  \fi \next}

\def\pto{.}

\def\multical #1{\def\arg{#1}%
  \ifx\arg\pto \let\next\relax
  \else
  \def\next{\expandafter
    \def\csname cal#1\endcsname{{\cal #1}}%
    \multical}%
  \fi \next}


\def\multimathop #1 {\def\arg{#1}%
  \ifx\arg\pto \let\next\relax
  \else
  \def\next{\expandafter
    \def\csname #1\endcsname{\mathop{\rm #1}\nolimits}%
    \multimathop}%
  \fi \next}

\multibold
qwertyuiopasdfghjklzxcvbnmQWERTYUIOPASDFGHJKLZXCVBNM.

\multical
QWERTYUIOPASDFGHJKLZXCVBNM.

\multimathop
dist div dom meas sign supp .


\def\accorpa #1#2{\eqref{#1}--\eqref{#2}}
\def\Accorpa #1#2 #3 {\gdef #1{\eqref{#2}--\eqref{#3}}%
  \wlog{}\wlog{\string #1 -> #2 - #3}\wlog{}}


\def\infess{\mathop{\rm inf\,ess}}
\def\supess{\mathop{\rm sup\,ess}}

\def\tonde #1{\left(#1\right)}

\def\graffe #1{\mathopen\{#1\mathclose\}}

\def\<#1>{\mathopen\langle #1\mathclose\rangle}
\def\norma #1{\mathopen \| #1\mathclose \|}

\def\iot {\int_0^t}
\def\ioT {\int_0^T}
\def\iO{\int_\Omega}
\def\intQt{\int_{Q_t}}
\def\intQ{\int_Q}

\def\dt{\partial_t}
\def\dn{\partial_\nu}
\def\ds{\,ds}

\def\cpto{\,\cdot\,}

\def\checkmmode #1{\relax\ifmmode\hbox{#1}\else{#1}\fi}
\def\aeO{\checkmmode{a.e.\ in~$\Omega$}}
\def\aeQ{\checkmmode{a.e.\ in~$Q$}}
\def\aeS{\checkmmode{a.e.\ on~$\Sigma$}}

\def\aat{\checkmmode{for a.a.~$t\in(0,T)$}}


\def\erre{{\mathbb{R}}}




\def\genspazio #1#2#3#4#5{#1^{#2}(#5,#4;#3)}
\def\spazio #1#2#3{\genspazio {#1}{#2}{#3}T0}

\def\L {\spazio L}
\def\H {\spazio H}
\def\W {\spazio W}

\def\C #1#2{C^{#1}([0,T];#2)}

\def\Vp{V^*}


\def\Lx #1{L^{#1}(\Omega)}
\def\Hx #1{H^{#1}(\Omega)}
\def\Wx #1{W^{#1}(\Omega)}

\def\Ldue{\Lx 2}
\def\Linfty{\Lx\infty}

\def\Huno{\Hx 1}
\def\Hdue{\Hx 2}


\def\LQ #1{L^{#1}(Q)}


\let\theta\vartheta
\let\eps\varepsilon

\let\TeXchi\chi                         
\newbox\chibox
\setbox0 \hbox{\mathsurround0pt $\TeXchi$}
\setbox\chibox \hbox{\raise\dp0 \box 0 }
\def\chi{\copy\chibox}


\def\muz{\mu_0}
\def\rhoz{\rho_0}

\def\muzsig{\mu_{0\sigma}}
\def\rhozsig{\rho_{0\sigma}}
\def\xizsig{\xi_{0\sigma}}

\def\rhomin{\rho_*}
\def\rhomax{\rho^*}
\def\ximin{\xi_*}
\def\ximax{\xi^*}
\def\rhosmrmp{(\rhos-\rhomax)^+}
\def\rhozsmrmp{(\rhozsig-\rhomax)^+}

\def\normaV #1{\norma{#1}_V}

\def\kamurho{\kappa(\mu,\rho)}

\def\kmin{\kappa_*}
\def\kmax{\kappa^*}

\def\coeff{1+2g(\rho)}

\def\mus{\mu_\sigma}
\def\rhos{\rho_\sigma}
\def\xis{\xi_\sigma}
\def\us{u_\sigma}
\def\kamurhosig{\kappa(\mus,\rhos)}
\def\kappat{\bar\kappa}
\def\Ku{K_1}
\def\Kd{K_2}

\def\pst{\pi^*}
\def\gst{g^*}



\Begin{document}


\title{\bf A vanishing diffusion limit 
in a nonstandard system of phase field equations%
}
\author{}
\date{}
\maketitle
\begin{center}
\vskip-2cm
{\large\bf Pierluigi Colli$^{(1)}$}\\
{\normalsize e-mail: {\tt pierluigi.colli@unipv.it}}\\[.2cm]
{\large\bf Gianni Gilardi$^{(1)}$}\\
{\normalsize e-mail: {\tt gianni.gilardi@unipv.it}}\\[.2cm]
{\large\bf Pavel Krej\v{c}\'{\i}$^{(2)}$}\\
{\normalsize e-mail: {\tt krejci@math.cas.cz}}\\[.2cm]
{\large\bf J\"urgen Sprekels$^{(3)}$}\\
{\normalsize e-mail: {\tt sprekels@wias-berlin.de}}\\[.4cm]
$^{(1)}$
{\small Dipartimento di Matematica ``F. Casorati'', Universit\`a di Pavia}\\
{\small via Ferrata 1, 27100 Pavia, Italy}\\[.2cm]
$^{(2)}$
{\small Institute of Mathematics, Czech Academy of Sciences}\\
{\small \v{Z}itn\'a 25, CZ-11567 Praha 1, Czech Republic}\\[.2cm]
$^{(3)}$
{\small Weierstra\ss-Institut f\"ur Angewandte Analysis und Stochastik}\\
{\small Mohrenstra\ss e\ 39, 10117 Berlin, Germany}\\[.8cm]
\end{center}

\Begin{abstract}
We are concerned with a nonstandard phase field model of Cahn-Hilliard type. 
The model, which was introduced {\gianni by Podio-Guidugli (Ric.\ Mat.\ 2006)}, 
describes two-species phase segregation 
and consists of a system of two highly nonlinearly coupled PDEs. 
It has been recently investigated 
{\gianni by Colli, Gilardi, Podio-Guidugli, and Sprekels in a series of papers:
see, in particular, SIAM J.\ Appl.\ Math.\ 2011
and Boll.\ Unione Mat.\ Ital.\ 2012}. 
{\GP In the latter contribution, the authors can treat 
the very general case in which the diffusivity coefficient of the parabolic PDE 
is allowed to depend nonlinearly on both variables.
In the same framework, this paper investigates the asymptotic limit 
of the solutions to the initial-boundary value problems
as the diffusion coefficient $\,\sigma\,$ 
in the equation governing the evolution of the order parameter tends to zero.
We prove that such a limit actually exists and solves the limit problem,
which couples a nonlinear PDE of parabolic type with an ODE accounting for the phase dynamics.
In the case of a constant diffusivity, we are able to show uniqueness 
and to improve the regularity of the solution.}
\\[2mm]
\End{abstract}


\salta

\pagestyle{myheadings}
\newcommand\testopari{\sc Colli \ --- \ Gilardi \ --- \ Krej\v{c}\'{\i} \ --- \ Sprekels}
\newcommand\testodispari{\sc Vanishing diffusion limit in nonstandard phase field systems}
\markboth{\testodispari}{\testopari}

\finqui


\section{Introduction}
\label{Intro}
\setcounter{equation}{0}
In this paper, we consider the following system 
\Bsist
  && \bigl( 1 + 2g(\rho) \bigr) \, \dt\mu
  + \mu \, g'(\rho) \, \dt\rho
  - \div \bigl( \kappa(\mu,\rho)\nabla\mu \bigr) = 0
  \label{Iprima}
  \\
  && \dt\rho - \sigma \Delta\rho + f'(\rho) = \mu \, g'(\rho)
  \label{Iseconda}
  \\
  && \big(\kappa(\mu,\rho)\nabla\mu\big)\cdot\nu|_\Gamma = 0
  \aand
  \dn\rho|_\Gamma = 0
  \label{Ibc}
  \\
  && \mu(0) = \muz
  \aand
  \rho(0) = \rhoz ,
  \label{Icauchy}
\Esist
\Accorpa\Ipbl Iprima Icauchy
in the unknown fields $\mu$ and~$\rho$, where the partial differential 
equations  \accorpa{Iprima}{Iseconda} {are meant {to hold} in a three-dimensional 
bounded domain $\Omega$, endowed with a smooth boundary~$\Gamma$, and in some time 
interval~$(0,T)$}. Relations \eqref{Icauchy} specify the initial conditions for 
$\mu$ and~$\rho$, while \eqref{Ibc} are nothing but homogeneous boundary 
conditions of Neumann type, involving precisely those boundary operators that match 
the elliptic differential operators in \accorpa{Iprima}{Iseconda}. 

This system has been recently addressed in the paper \cite{CGPS7}: the 
existence of solutions has been proved, thus complementing and
extending the results of the papers {\GP\cite{CGPS3, CGPS4, CGPS6}} concerned with simpler or reduced versions of the problem. 

Here, we are interested to investigate the asymptotic \bhv\ of the above initial-boundary value problem \Ipbl\ 
as the positive diffusion coefficient $\sigma$ appearing in \eqref{Iseconda} tends to~$0$.

Let us briefly explain the modelling background for \Ipbl. Such a system comes
from a generalization of the phase-field model of viscous Cahn-Hilliard type originally proposed in \cite{Podio}, 
and it aims to describe the phase segregation of two species 
(atoms and vacancies, say) 
on~a lattice in presence of diffusion. 
The state variables are the  {\sl order parameter\/} $\rho$,
interpreted as {the volume density of one of the two species},  
and the {\sl chemical potential\/}~$\mu$. 
For physical reasons, $\mu$~is required to be nonnegative, 
while the phase parameter $\rho$ must of course take values in the domain of $f'$. 

We also recall the features of \cite{CGPS3} and what has been generalized in \cite{CGPS6, CGPS7}. 
Firstly, the nonlinearity $f$ considered in \cite{CGPS3} is 
a double-well potential defined in~$(0,1)$, 
whose derivative $f'$ {diverges} at the endpoints $\rho=0$ and $\rho=1$: e.g., for
$f=f_1+f_2$ with $f_2$ smooth, 
one can take  
\Beq
  f_1(\rho)=c\,(\rho\,\log(\rho)+(1-\rho)\,\log(1-\rho)),  
  \label{ln} 
\Eeq
with $c$ a positive constant. 
In this paper, we let $f_1 :\erre \to [0,+\infty] $ be a convex, proper and lower semicontinuous function 
so that its subdifferential (and not the derivative) 
is a maximal monotone graph from $\erre$ to~$\erre$. 
Then, we rewrite equation~\eqref{Iseconda} as a differential inclusion, 
in which the derivative of the convex part $f_1$ of $f$ is replaced by 
the subdifferential $\beta:=\partial f_1$,~i.e.,
\Beq
  \dt\rho - \sigma \Delta\rho + \xi + f_2'(\rho) = \mu g'(\rho)
  \quad \hbox{with} \quad
  \xi \in \beta(\rho).
  \label{Isecondabis}
\Eeq 
Note that $f_1$ need not be differentiable in its domain, 
so that its possibly nonsmooth and multivalued subdifferential  $\beta:=\partial f_1$ 
appears in \eqref{Iseconda} in place of~$f'_1$. 
In general, $\beta$ is only a graph, 
not necessarily a function, and it may include vertical 
(and horizontal) 
lines, as for example when
\begin{equation}
  \label{ex2}
  f_1 (\rho) =  I_{[0,1]} (\rho) =  \ \left\{
  \begin{array}{ll}
    0
    & \text{if \ $0\leq \rho \leq 1$}
    \\[0.1cm]
    +\infty \ 
    & \text{elsewhere}
  \end{array}
  \right. 
\end{equation}
and $\beta = \partial I_{[0,1]}$ is specified by
\begin{equation}
  \label{ex0}
  \xi \in \beta (\rho) 
  \quad \hbox{ if and only if } \quad 
  \xi \ \left\{
  \begin{array}{ll}
    \displaystyle
    \leq \, 0 \
    &\hbox{if } \ \rho=0   
    \\[0.1cm]
    = \, 0 \
    &\hbox{if } \ 0< \rho < 1  
    \\[0.1cm]
    \geq \, 0 \
    &\hbox{if } \  \rho  = 1  
    \\[0.1cm]
  \end{array}
  \right. .
\end{equation}
Secondly, while in \cite{CGPS3} $g$ was simply taken as the identity map $g(\rho)=\rho$, 
in \cite{CGPS6,CGPS7} $g$ is allowed be any nonnegative smooth function, 
defined (at~least) in the domain where $f_1$ and its subdifferential live. 
The presence of such a function $g$ allows 
for a more general \bhv\ of (the related term in) the free energy, 
which reads 
\begin{equation}
  \label{fe-2}
  \psi(\rho,\nabla\rho,\mu) 
  = - \frac\mu2 - \mu \, g(\rho) + f (\rho) + \frac{\sigma}{2}|\nabla\rho|^2.
\end{equation}
Indeed, in particular $g(\rho)$ is not obliged, as it was instead for $g(\rho)=\rho$, 
to take its minimum value at $\rho=0$, be increasing and with maximum value at $\rho=1$ 
(when $D(f_1)= [0,1]$), but we may have many other instances like, e.g., 
a specular \bhv\ of $g$ around the extremal points of the domain~of~$f$. 
Here, we have to impose an additional restriction on~$g$,
which however looks reasonable from the modelling {\gianni point of view}: 
we postulate that $g$ is~a (smooth) concave function, 
which in turn implies convexity with respect to $\rho$ of the term  $-\mu \, g(\rho)$ in the free energy \eqref{fe-2}. 
However, let us recall that $ f $ may stand for a multi-well 
potential in which the nonconvex perturbations are incorporated into $f_2$, 
so that $\psi$ in its entirety needs not be convex with respect to $\rho$.

An important generalization that is considered in this paper concerns the diffusivity~$\kappa$. 
In~\cite{CGPS3}, $\kappa$~was just assumed to be a constant function, 
but it can be a positive-valued, continuous, bounded, and nonlinear function of~$\mu$ 
(and this was the setting of~\cite{CGPS6}), 
or of $\mu$ and $\rho$ as it is postulated in~\cite{CGPS7}. 
For simplicity, 
we confine ourselves to study of the convergence properties of the solution 
under an assumption that guarantees uniform parabolicity, i.e., $\kappa\geq\kmin>0$. 
We point out that \cite{CGPS6} treats the situation of  $\kappa$ 
depending only on $\mu$ and possibly degenerating somewhere.

Therefore, the system 
\Bsist
  && \bigl( 1 + 2g(\rho) \bigr) \, \dt\mu
  + \mu \, g'(\rho) \, \dt\rho
  - \div \bigl( \kappa(\mu,\rho)\nabla\mu \bigr) = 0
  \label{Iprimater}
  \\
  &&
  \dt\rho - \sigma \Delta\rho + \xi + f_2'(\rho) = \mu g'(\rho)
  \quad \hbox{with} \quad
  \xi \in \beta(\rho), 
  \label{Isecondater}
  \\
  && \big(\kappa(\mu,\rho)\nabla\mu\big)\cdot\nu|_\Gamma = 0
  \aand
  \dn\rho|_\Gamma = 0
  \label{Ibcter}
  \\
  && \mu(0) = \muz
  \aand
  \rho(0) = \rhoz ,
  \label{Icauchyter}
\Esist
turns out the initial and boundary value problem for a nonstandard 
and highly nonlinear phase field system in which however 
the role usually played by the temperature is here conducted by the chemical potential~$\mu$. 
In the study of  phase field systems, it has been always considered rather important 
to analyze the \bhv\ of the problem as the coefficient $\sigma $ of the diffusion term 
in the phase parameter equation tends to~$0$. 
The limiting case $\sigma=0$ corresponds indeed to a pointwise ordinary differential equation (or inclusion) 
\Beq
  \dt\rho + \xi + f_2'(\rho) = \mu g'(\rho), 
  \quad 
  \xi \in \beta(\rho), 
  \label{Isecondaquater}
\Eeq
in place of \eqref{Isecondater}, and to an expression for the free energy \eqref{fe-2} 
in which the last quadratic term accounting for {\gianni nonlocal} interactions is removed. 

In fact, especially for the choice \accorpa{ex2}{ex0}, the limiting problem
can be formulated in terms of hysteresis operators: 
in particular, the so-called {\it stop} and {\it play} operators are involved; 
the interested reader can find some discussion and various results 
on this class of problems in \cite{CS1, CS2, GKS, KS1, KS2, KS3, KSZ}.

By collecting a number of estimates independent of $\sigma$ for the solution 
$(\mus, \rhos)$ to the problem \accorpa{Iprimater}{Icauchyter}, 
by weak and weak star compactness we {\gianni prove} that any limit 
in a suitable topology of a convergent subsequence of $\{ (\mus, \rhos)\}$
yields a solution to the limiting problem in which \eqref{Isecondater}
is replaced by \eqref{Isecondaquater}. 
{\GP
Furthermore, under natural compatibility conditions on the nonlinearities and the initial data,
we show boundedness for all the components of any solution to the limit problem.
Finally, in the special case of a constant mobility~$\kappa$ in~\eqref{Iprimater},
we prove that the solution is unique and more regular than required.}

The paper is \organiz ed as follows. 
In the next section, we state precise assumptions along with our results.
The basic a priori estimates independent of $\sigma$ are proved in Section~\ref{Est-lim} 
and they allow us to pass to the limit 
by compactness and monotonicity techniques. 
{\GP
Finally, the last section is devoted to the study of the limit problem
and our boundedness, uniqueness, and further regularity properties are proved.}


\section{Assumptions and results}
\label{MainResults}
\setcounter{equation}{0}

The aim of this section is to introduce precise assumptions on the data 
for the mathematical problem under investigation, and establish our main result. 
We assume $\Omega$ to be a bounded connected 
open set in $\erre^3$ with smooth boundary~$\Gamma$
({treating} lower-dimensional cases would require {only }minor changes) 
and let  $T\in(0,+\infty)$ stand for a final time. 
We introduce the spaces
\Beq
  V := \Huno,
  \quad H := \Ldue , \quad
  W := \graffe{v\in\Hdue:\ \dn v = 0 \ \hbox{on $\Gamma$}}
  \label{defspazi}
\Eeq
and endow them with their standard norms,
for which we use a self-explanato\-ry notation like $\normaV\cpto$. 
For powers of these spaces, norms are denoted by the same symbols.
We remark that the embeddings $W\subset V\subset H$ are compact,
because $\Omega$ is bounded and smooth. 
The symbol $\<\cpto,\cpto>$ denotes the duality product 
between~$\Vp$, the dual space of~$V$, and~$V$ itself.
Moreover, for $p\in[1,+\infty]$, we write $\norma\cpto_p$ for the usual norm in~$L^p(\Omega)$; 
as no confusion can arise, the symbol $\norma\cpto_p$ is used for the norm in $L^p(Q)$ as well,
where $Q:=\Omega\times(0,T)$.

Now, we present the structural assumptions we make. 
It is useful to fix an upper bound for~$\sigma$, that~is,
\Beq 
  \label{siguno}
  0<\sigma \leq 1 .
\Eeq
Then, for the diffusivity coefficient $\kappa$ we assume that
\Bsist
  \hskip-1.5cm
  & \hbox{$\kappa:(m,r)\mapsto\kappa(m,r)$ is continuous from $[0,+\infty)\times\erre$ to $\erre$}{,}
  &
  \label{hpk}
  \\
  \hskip-1.5cm& \hbox{the partial derivatives $\partial_r\kappa$ and $\partial_r^2\kappa$
  exist and are continuous}{,}
  &
  \label{hpdk}
  \\
  \hskip-1.5cm& \kmin,\, \kmax \in (0,+\infty){,}
  &
  \label{hpcost}
  \\
  \hskip-1.5cm& \displaystyle \kmin \leq \kappa(m,r) \leq \kmax, \enskip
  |\partial_r\kappa(m,r)| \leq \kmax,
  \enskip 
  |\partial_r^2\kappa(m,r)| \leq \kmax \quad 
  \hbox{for $m\geq0$ and $r\in\erre$,}
  &
  \label{hpkbis}
\Esist
and for other nonlinearities we require that
\Bsist
  \hskip-1.5cm&
  f = f_1 + f_2 \,, \quad f_1:\erre \to [0,+\infty], \quad f_2 :\erre \to \erre, 
  &
  \label{hpfg}
  \\
  \hskip-1.5cm& \hbox{$f_1$ is convex, proper, l.s.c. and $f_2$ is a $C^2$ function},&
  \label{hpf12}
   \\
  \hskip-1.5cm& \hbox{$ g\in C^2 (\erre)$, \quad $g(r)\geq0$ \  and \ $g''(r) \leq 0$ \ for  $r\in\erre$,}&
  \label{hpg}
  \\
  \hskip-1.5cm& \hbox{$f_2'$, $g$, and $g'$ are Lipschitz continuous}. & 
  \label{hpfdueg}
\Esist
It is convenient to introduce the notations
\Bsist
  \hskip-1.5cm& \kappa' := \partial_r \kappa , \quad
  \kappa'' := \partial_r^2 \kappa , \quad
  \beta := \partial f_1 \,,
  \aand
  \pi := f_2' &
  \label{defbp}
  \\
  \hskip-1.5cm& \displaystyle K(m,r) := \int_0^m \!\!\! \kappa(s,r) \ds, \
  \Ku(m,r) := \int_0^m \!\!\! \kappa'(s,r) \ds, \
  \Kd(m,r) := \int_0^m \!\!\! \kappa''(s,r) \ds &
  \non
  \\
  \hskip-1.5cm& \hskip9cm \hbox{for $m\geq0$ and $r\in\erre$}.&
  \label{defK}
\Esist
\Accorpa\Hpstruttura hpk defK
We denote by $D(f_1)$ and $D(\beta)$ the effective domains
of~$f_1$ and~$\beta$, respectively. Thanks to \eqref{hpkbis}, it is clear~that
\Beq
  \max \{ |K(m,r)| , |\Ku(m,r)| , |\Kd(m,r)| \}
  \leq \kmax m
  \quad \hbox{for every $m\geq0$ and $r\in\erre $}.
  \label{lingr}
\Eeq
We also note that the structural assumptions of~\cite{CGPS6}
are fulfilled if $\kappa$ only depends on~$m$, 
and that, due to the presence of $\beta(\rho)$, 
a~strong singularity {\gianni in equation \eqref{Isecondater}} is allowed. 
On the other hand,  equation \eqref{Iprimater} is uniformly parabolic, 
since $g$ is nonnegative and $\kappa$ is bounded away from zero.

\Brem
\label{Justif}
Let us recall that any convex, proper, l.s.c.\ function is bounded from 
below by an affine function (cf., e.g., \cite[Prop.~2.1, p.~51]{Barbu}),
whence the assumption $f_1 \geq 0 $ looks reasonable, as one can suitably modify the smooth perturbation~$f_2$. 
Moreover, we point out that the sign conditions $g\geq 0$ and $g''\leq 0$ 
are just needed on the set~$D(\beta)$, for $g$ can be extended outside of $D(\beta)$ accordingly. 
\Erem

Concerning the initial data, we~require~that
\Bsist
  & \muz \in V , \quad
  \muz \geq 0  \quad \aeO,&
  \label{hpmuzero}
  \\
  & \rhoz \in V , \quad
\rhoz \in D(f_1) \quad \aeO , \quad f_1 (\rhoz ) \in L^1(\Omega)&
  \label{hprhozero}
\Esist
\Accorpa\Hpdati hpmuzero hprhozero
and point out that the above assumptions regard the initial data for the limiting problem, 
i.e., the one with \eqref{Isecondaquater} in place of~\eqref{Isecondater}. 
On the other hand, let us consider a family of initial data $\muzsig, \, \rhozsig $ with
\Bsist
  \hskip-1.5cm& \muzsig \in V \cap \Linfty , \quad
  \muzsig \geq 0
  \quad \aeO,&
  \label{hpmuzsig}
  \\
  \hskip-1.5cm&\rhozsig \in W ,
  \quad \hbox{there is $\ \xizsig \in H \ $ such that } \
  \rhozsig \in D(\beta), \ \,
  \xizsig \in\beta(\rhozsig) \ \ \aeO,
  &
  \label{hprhozsig}
\Esist
\Accorpa\Hpdatisig hpmuzsig hprhozsig
that approximate $\muz, \, \rhoz$ in the sense that  
\Bsist
  & \muzsig \to \muz \ \hbox{ and }\ \rhozsig \to \rhoz \ \hbox{ weakly in } \ V,
  &
  \label{P1}
  \\
  \quad
  & \norma{f_1 (\rhozsig )}_1 \ \hbox{ is bounded independently of } \ \sigma.
  &
  \label{P2}
\Esist
For the reader's convenience, we show that such a family  $\{ \muzsig, \, \rhozsig \} $ actually exists. 
Of course, if $\muz \not\in \Linfty$ we can take as 
$ \muzsig $ some truncation of $\muz$, e.g., $ \muzsig = \min \{ \muz , 1/\sigma \}$.
Concerning $ \rhozsig $, one possible choice is letting $ \rhozsig \in W $ denote the solution~to 
\Beq
 \rhozsig  - \sigma \Delta\rhozsig + \sigma \xizsig = \rhoz , 
  \quad \hbox{with} \quad
  \xizsig \in \beta(\rhozsig),  \ \ \aeO. 
  \label{P3}
\Eeq
Indeed, the elliptic problem \eqref{P3} has a unique solution for all $\sigma >0$, 
since $-\Delta+\beta $ is a maximal monotone graph in $H\times H$ with effective domain 
$$
  \{ v\in W : \  \exists\  \eta \in H \hbox{ such that } v \in D(\beta), \ \eta \in \beta(v) \  \aeO \}.
$$
Thus, $\rhozsig$ is nothing but the outcome of the application of the resolvent of $-\Delta+\beta$ to~$\rhoz$ 
(let~us refer to \cite{Barbu} and \cite{Brezis} for basic definitions and properties of maximal monotone operators). 
A~formal test of the equality in \eqref{P3} by $\xizsig$ and the definition of subdifferential 
lead us to the estimate 
\Beq
  \int_\Omega f_1( \rhozsig) + \sigma \norma{ \xizsig}^2_H \leq  \int_\Omega f_1( \rhoz) ,  
  \label{P4}
\Eeq 
which ensures \eqref{hprhozsig} and \eqref{P2}, thanks to the {\gianni nonnegativity} of~$f_1$. 
A~rigorous way of proving the existence of $\rhozsig$ and estimate~\eqref{P4} 
passes through the use of the Yosida approximation~$\beta_\sigma$ 
(see, e.g., \cite[p.~28]{Brezis}) 
in place of~$\beta$.   

Now, we recall the result proved in \cite{CGPS7} that allows us to specify a solution 
to the problem~\eqref{Iprimater}--\eqref{Ibcter}, with $\sigma>0$, 
which fulfills the appropriate initial conditions. 

\Bprop
\label{Esisigma}
Assume that both \Hpstruttura\ and \Hpdatisig\ hold.
Then, there exists at least one triplet $(\mus,\rhos,\xis)$ satisfying 
\Bsist
  \hskip-1.5cm& \rhos \in \W{1,\infty}H \cap \H1V \cap \L\infty W
  \label{regrho}, &
  \\
  \hskip-1.5cm & 
  \xis \in \L\infty H  ,&
  \label{regxi}
  \\
  \hskip-1.5cm& \mus \in \H1H \cap \L\infty V \cap \LQ\infty , \quad
  \mus \geq 0 \quad \aeQ,&
  \label{regmu}
  \\
  \hskip-1.5cm& \div \bigl( \kamurhosig\nabla\mus \bigr) \in L^2(Q) 
  \aand
  \bigl( \kamurhosig\nabla\mu \bigr) \cdot \nu = 0
  \quad \aeS,&
  \label{regdiv}
\Esist
and solving the system of equations and conditions in the following strong form
\Bsist
  \hskip-1.5cm & \bigl( 1+ 2 g(\rhos) \bigr) \dt\mus + \mus \, g'(\rhos) \, \dt\rhos
  - \div \bigl( \kamurhosig \nabla\mus \bigr) = 0
  & \quad \aeQ,
  \label{prima}
  \\
  \hskip-1.5cm & \dt\rhos - \sigma \Delta\rhos + \xis + \pi(\rhos)
  = \mus \, g'(\rhos)
  \aand \xis \in \beta(\rhos)
  & \quad \aeQ,
  \label{seconda}
  \\
  \hskip-1.5cm & \mus(0) = \muzsig
  \aand
  \rhos(0) = \rhozsig
  & \quad \aeO.
  \label{cauchy}
\Esist
\Accorpa\Pbl prima cauchy
\Eprop

Let us point out that equation \eqref{prima} can be rewritten as 
\Bsist
  \dt \us - \div \bigl( \kamurhosig \nabla\mus \bigr) = \mus \, g'(\rhos) \, \dt\rhos , \qquad
  \non
  \\
  \quad \hbox{where} \quad 
  \us = ( 1+ 2 g(\rhos)) \mus,
  \quad \aeQ,
  \label{var-prima}
\Esist
and the auxiliary variable $\us$ has been added. 
Now, we take advantage of a variational formulation of \eqref{var-prima}
which also accounts for the boundary condition in \eqref{regdiv}, that~is,
\Bsist
  \< \dt \us (t), v > +\int_\Omega  \bigl( \kamurhosig \nabla\mus \bigr)(t) \cdot \nabla v
  = \int_\Omega \mus \, g'(\rhos) \, \dt\rhos \, v \qquad
  \non
  \\
  \quad \hbox{for all } v\in V \hbox{ and a.a. } t\in (0,T).
  \label{pri-variaz}
\Esist
{\gianni The main result of this paper reads as follows.}

\Bthm
\label{Convergenza}
Assume that~\Hpstruttura\ and~\eqref{hpmuzero}--\eqref{P2} hold. 
For any $\sigma \in (0,1]$ let $(\mus,\rhos,\xis)$ be the 
triplet defined by Proposition~\ref{Esisigma} and let 
$\us : = (1+ 2 g(\rhos)) \mus$. 
Then, there exists a subsequence, still labelled by the parameter~$\sigma$, 
and a quadruplet  $(\mu,\rho,\xi, u)$ such~that
\Bsist 
\hskip-1cm& \mus \to \mu
  & \quad \hbox{weakly star in $\L\infty H\cap\L2V$}, 
  \qquad
  \label{convmu}
  \\
  \hskip-1cm& \rhos \to \rho
  & \quad \hbox{weakly star in $\H1H\cap\L\infty V$},
  \label{convrho}
  \\
  \hskip-1cm& \xis \to \xi
  & \quad \hbox{weakly in $ L^2(Q)$},
  \label{convxi}
  \\
  \hskip-1cm& \us \to u  
  & \quad \hbox{weakly in $W^{1, 4/3} (0,T;\Vp)\cap\L2{\Wx{1,3/2}}$}
  \label{convu}
\Esist
as $\sigma \searrow 0.$ Moreover, any quadruplet  $(\mu,\rho,\xi, u)$ 
that is found as limit of converging subsequences yields  
a solution to the following limit problem 
\Bsist
  & \< \dt u(t), v > + \int_\Omega \blu{\kamurho \nabla\mu} (t) \cdot \nabla v
  = \int_\Omega \mu \, g'(\rho) \, \dt\rho \, v \qquad& \non\\
  &\hskip6cm \hbox{for all } v\in V \hbox{ and a.a. } t\in (0,T),&
  \label{primau}
  \\
   & u= (1+2g(\rho))\mu 
  \quad \aeQ,&\label{defu}
  \\
   &   \dt\rho + \xi + \pi(\rho)
   = \mu \, g'(\rho)
  \aand \xi \in \beta(\rho)
  \quad \aeQ, & \label{secondez}
  \\
  & \mu(0) = \muz
  \aand
  \rho(0) = \rhoz
   \quad \aeO.& 
  \label{cauchyz}
\Esist
\Accorpa\Pblz primau cauchyz
\Ethm

\Brem
\label{Sign}
The nonnegativity property $\mu \geq 0 $ \aeQ\ plainly follows from \eqref{regmu} and~\eqref{convmu}. 
\Erem

\Brem
\label{Stop}
One standard situation for the limit problem~\Pblz\
is obtained for $\beta = \partial I_{[0,1]}$ (cf.~\eqref{ex2}--\eqref{ex0}).
In this case \eqref{secondez} becomes 
\Beq
   - \pi(\rho)
   + \mu \, g'(\rho) - \dt\rho  \in \partial I_{[0,1]}(\rho)
  \quad \aeQ.  \label{P9}
\Eeq
Then, if one introduces the generalized ``freezing index'' 
$$
  w (x,t) := \int_0^t  (- \pi(\rho)
  + \mu \, g'(\rho)) (x,s) ds, \quad (x,t)\in Q, 
$$
we thus have $\dt w - \dt\rho \in  \partial I_{[0,1]}(\rho)$, or equivalently, 
$\rho = {\cal S}_K [w]$, 
where  $ {\cal S}_K $ is the stop hysteresis operator associated 
with the closed convex set $K=[0,1]$ (see, e.g., \cite{KS1, KS2, KS3}). 
Hence, we may rewrite \eqref{P9}~as 
\Beq
 \dt w  =  - \pi({\cal S}_K [w])
   +  \mu \, g'({\cal S}_K [w]) 
  \quad \aeQ.  \non
\Eeq
\Erem

In addition to the convergence result stated in Theorem~\ref{Convergenza},
one {\gianni can derive} boundedness for both the components $\rho$ and $\xi$
of any solution to the limit problem,
provided that special additional requirements are satisfied,
namely, by assuming that there exist real constants 
$\rhomin,\,\rhomax,\,\ximin,\,\ximax$ such~that  
\Bsist
  & \rhomin,\,\rhomax \in D(\beta) , \quad
  \ximin \in \beta (\rhomin) , \quad
  \ximax \in \beta (\rhomax),
  & \label{P13}
  \\
  & \ximin + \pi ( \rhomin ) \leq 0 , \quad \ximax + \pi ( \rhomax ) \geq 0 , 
  & \label{P14}
  \\
  & g'(\rhomin) \geq 0 , \quad g'(\rhomax) \leq 0.
  & \label{P15}
\Esist
\Accorpa\Perlimitatezza {P13} {P15}

\Bthm
\label{Linftyz} 
{\GP In addition to the assumptions of Theorem~\ref{Convergenza},
suppose that \Perlimitatezza\ and}
\Beq
  \rhomin \leq \rhoz \leq  \rhomax \quad \aeO   
  \label{eq:PJ1}
\Eeq
{\GP hold}.
Then, the components $\rho$ and $\xi$ of any solution 
$(\mu,\rho,\xi, u)$ to problem~\Pblz\ satisfy 
\Beq
 \rhomin \leq \rho \leq  \rhomax
 \aand
 \ximin \leq \xi \leq \ximax
 \quad \aeQ .
 \label{linftyz}
\Eeq
{\GP If moreover 
\Beq
  \muz \in \Linfty 
  \label{hpmuzbdd}
\Eeq
and $\kappa = \kappa_0$ is constant,
then the solution of Problem \eqref{primau}--\eqref{cauchyz} is unique~and 
\Beq
  {\GP \mu \in \H1H \cap \L\infty V \cap \L2W .}
  \label{piuregolare}
\Eeq
}%
\Ethm

\Brem
\label{Maxpri}
We observe that the above result is very general.
Indeed, assumptions \Perlimitatezza\ are fulfilled with suitable constants
for any graph~$\beta$ with bounded domain that \generaliz es the examples \eqref{ln} or~\eqref{ex2}.
Of course, the decreasing function~$g'$ (cf.~\eqref{hpg}) should not assume a definite sign on~$D(\beta)$.
\Erem

Now, we list a number of tools and notations we owe to throughout the paper.
We repeatedly use the elementary Young inequalities
\Bsist
  && a\,b\leq \gamma a^2 + \frac 1 {4\,\gamma}\,b^2
  \aand
  a\,b \leq \theta a^{\frac 1\theta} + (1-\theta) b^{\frac 1{1-\theta}}
  \non
  \\
  && \quad \hbox{for every $a,b\geq 0$, $\gamma>0$, and $\theta\in(0,1)$}
  \label{young}
\Esist
as well as the H\"older and Sobolev inequalities.
The precise form of the latter we use is the following
\Bsist
  \hskip-1cm & \Wx{1,p} \subset \Lx q
  \aand
  \norma v_{q} \leq C_{p,q} \norma v_{\Wx{1,p}}
  \quad \hbox{for every $v\in\Wx{1,p}$,} &
  \non
  \\
  \hskip-1cm & \quad \hbox{provided that} \quad
  1 \leq p <3
  \aand \displaystyle
  1 \leq q \leq p^* := \frac {3p} {3-p} &
  \label{sobexp}
\Esist
with a constant $C_{p,q}$ in \eqref{sobexp} depending only on $\Omega$, $p$, and~$q$,
since $\Omega\subset\erre^3$.
Moreover
\Beq
  \hbox{the embedding} \quad
  \Wx{1,p} \subset \Lx q
  \quad \hbox{is compact if} \quad 
  1 \leq q < p^* .
  \label{compsobolev}
\Eeq
The particular case $p=2$ of \eqref{sobexp} becomes
\Beq
  V \subset \Lx q
  \aand
  \norma v_{q} \leq C \normaV v
  \quad \hbox{for every $v\in V$ and $q\in[1,6]$}
  \label{sobolev}
\Eeq
where $C$~depends only on~$\Omega$.
Moreover, the compactness inequality
\Beq
  \norma v_q \leq \eps \norma{\nabla v}_2 + C_{q,\eps} \norma v_2
  \quad \hbox{for every $v\in V$, $q\in[1,6)$, and $\eps>0$}
  \label{compact}
\Eeq
holds for some constant $C_{q,\eps}$ depending on $\Omega$, $q$, and~$\eps$, only.
We also recall the interpolation inequalities, which hold for any $\theta\in[0,1]$,
\Bsist 
  && \norma v_r
  \leq \norma v_p^\theta \, \norma v_q^{1-\theta}
  \,\quad\forall \, v \in \Lx p \cap \Lx q,
  \non
  \\
  && \quad \hbox{where }\,p,q,r\in [1,+\infty]
  \aand 
  \frac 1 r =\frac \theta p + \frac {1-\theta} q \,.\qquad 
  \label{interpolx}
  \\
  && \norma v_{\L{r_1}{\Lx{r_2}}}
  \leq \norma v_{\L{p_1}{\Lx{p_2}}}^\theta \, \norma v_{\L{q_1}{\Lx{q_2}}}^{1-\theta}
  \vphantom \int
  \non
  \\
  && \quad\forall \, v \in \L{p_1}{\Lx{p_2}} \cap \L{q_1}{\Lx{q_2}},
  \non
  \\
  && \quad \hbox{where} \quad p_i\,,q_i\,,r_i\in [1,+\infty]
  \aand
  \frac 1 {r_i} =\frac \theta {p_i} + \frac {1-\theta} {q_i} \quad
  \quad \hbox{for $i=1,2$}.
  \qquad
  \label{interpolxt}
\Esist
We observe that \eqref{interpolx} implies 
$\norma v_r \leq \theta \norma v_p + (1-\theta) \norma v_q$
for every $v\in\Lx p\cap\Lx q$
thanks to the Young inequality,
and a similar remark holds for~\eqref{interpolxt}.
Thus, we have the continuous embeddings
\Beq
  \Lx p \cap \Lx q \subset \Lx r
  \aand
  \L{p_1}{\Lx{p_2}} \cap \L{q_1}{\Lx{q_2}} \subset \L{r_1}{\Lx{r_2}} .
  \non
\Eeq
We stress the important case of the space
$\L\infty\Ldue\cap\L2{\Lx6}$, which occurs several times in the sequel
and corresponds to
$p_1=\infty$, $p_2=2$, $q_1=2$, and $q_2=6$.
In particular, the choices $\theta=2/5$ and $\theta=1/7$ 
yield the inequalities (for~every $v$ of the above space) 
and~the continuous embeddings
\Bsist
  && \norma v_{\LQ{10/3}}
  \leq \norma v_X^{2/5} \norma v_Y^{3/5} 
  \aand
  X \cap Y \subset \LQ{10/3}
  \label{inter10-3-Q}
  \\
  && \norma v_{\L{7/3}{\Lx{14/3}}}
  \leq \norma v_X^{1/7} \norma v_Y^{6/7}
  \aand
  X \cap Y \subset \L{7/3}{\Lx{14/3}}
  \qquad
  \label{inter7-14-3}
  \\
  && \quad \hbox{where} \quad \vphantom\int
  X := \L\infty\Ldue 
  \aand
  Y := \L2{\Lx6} .
  \non
\Esist
Notice that we can take $v\in\L\infty H\cap\L2V$
in \accorpa{inter10-3-Q}{inter7-14-3},
since $V\subset\Lx6$.
Finally, we set 
\Beq
  Q_t := \Omega \times (0,t)
  \quad \hbox{for $t\in[0,T]$},
  \label{defQt}
\Eeq
and, again throughout the paper,
we use a small-case italic $c$ for different constants, that
may only depend 
on~$\Omega$, the final time~$T$, the shape of the nonlinearities $f$ and~$g$, 
and the properties of the data involved in the statements at hand; 
a~notation like~$c_\eps$ signals a constant that depends also on the parameter~$\eps$. 
The reader should keep in mind that the meaning of $c$ and $c_\eps$ might
change from line to line and even in the same chain of inequalities, 
whereas those constants we need to refer to are always denoted by 
capital letters, just like $C$ in~\eqref{sobolev}.


\section{The asymptotic analysis}
\label{Est-lim}
\setcounter{equation}{0}

In this section, we prove Theorem~\ref{Convergenza}, 
which ensures the existence of a solution to problem~\Pblz\ 
along with the convergence properties stated in~\eqref{convmu}--\eqref{convu}. 

Then, for any $\sigma \in (0,1]$ we let $(\mus,\rhos,\xis)$ denote 
the triplet defined by Proposition~\ref{Esisigma} and set 
$\us : = ( 1+ 2 g(\rhos)) \mus$. 
The existence of $(\mus,\rhos,\xis)$ has been proved in \cite{CGPS7}: 
we follow in parts the arguments developed there in order to recover 
useful estimates independent of~$\sigma$.
Before that, let us remark that the property $\mus \geq 0$ can be verified 
by simply multiplying  equation~\eqref{prima} by $-\mus^-$,
the negative part of~$\mus$, and integrate over~$Q_t$. 
In principle, in this computation one has to define $\kappa$ everywhere, 
e.g., by taking an even extension $\kappat$ with respect to the first variable. 
We observe that
\Beq
  \bigl[ \bigl( 1+2 g(\rhos (t)) \bigr) \, \dt\mus + \mus \, g(\rhos) \, \dt\rhos \bigr] \, (-\mus^-)
  = \frac 12 \, \dt \bigl( (1+2 g(\rhos (t))) \, |\mus^-|^2 \bigr).
  \non
\Eeq
Hence, by using $\muzsig\geq0$ 
and owing to the boundary condition in \eqref{regdiv}, we have
\Bsist
  \frac 12 \iO (1+2 g(\rhos (t))) \, |\mus^-(t)|^2
  + \intQt \kappat(\mus,\rhos) |\nabla\mus^-|^2 
  = 0 
  \quad \hbox{for a.a. $t\in (0,T)$}.
  \non
\Esist
As both $g$ and $\kappat $ are nonnegative,
this implies $\mus^-=0$, that is, $\mus\geq0$ a.e.\ in~$Q$.

\step 
First a priori estimate

We test \eqref{prima} by $\mus$ and point out that
\Beq
  \bigl[ \bigl( 1+2g(\rhos)\bigr) \, \dt\mus
  + \mus \, g'(\rhos) \, \dt\rhos \bigr] \mus
  = \frac 12 \, \dt \bigl[ (1+2g(\rhos) \mus^2 \bigr].
  \label{per1stima}
\Eeq
Thus, by integrating over~$(0,t)$, where $t\in[0,T]$ is arbitrary,
we obtain
\Beq
  \iO \bigl( 1+2g(\rhos(t)) \bigr) |\mus(t)|^2 + 2 \intQt \kappa(\mus(s),\rhos(s)) |\nabla\mus|^2
  = \iO (1+2g(\rhozsig)) \muzsig^2 \,.
  \non
\Eeq
We recall that $g$ is nonnegative and Lipschitz continuous (cf.~\accorpa{hpg}{hpfdueg}). 
Moreover, $\rhozsig, \,\muzsig$ are both uniformly bounded in $V$ by \eqref{P1}, whence
$$
  \iO (1+2g(\rhozsig)) \muzsig^2  \leq c \left( \norma{\muzsig}_2^2 + 
  \norma{\rhozsig}_2 \norma{\muzsig}_4^2 \right) \leq c 
$$
owing to the H\"older and Sobolev inequalities (see~\eqref{sobolev}). 
Then, in view of $g\geq0$ and $\kappa\geq\kmin>0$, 
from \eqref{per1stima} it follows that
\Beq
  \norma\mus_{\L\infty H} + \norma{\mus}_{\L2V} \leq c .
  \label{primastima}
\Eeq

\step 
Second a priori estimate

We add $\rhos$ to both sides of \eqref{seconda} and test by~$\dt\rhos$.
On account of \eqref{hpfg}--\eqref{hpf12} and \eqref{defbp},
we obtain
\Bsist
  && \intQt |\dt\rhos|^2
  + \frac 12 \, \norma{\rhos(t)}_H^2 + \frac {\sigma}2 \, \norma{\nabla \rhos(t)}_H^2
  + \iO f_1(\rhos(t))
  \non
  \\
  && = \frac {\sigma}2 \iO |\nabla\rhozsig|^2
  + \iO f (\rhozsig)
  + \frac 12 \iO \tonde{\rhos^2 (t) - 2 f_2(\rhos(t))} 
  + \intQt \mus g'(\rhos) \dt\rhos 
  \non
\Esist
for every $t\in[0,T]$. 
Then, thanks to the Lipschitz continuity of $f'_2$
and $g$, and owing to the bounds entailed by \accorpa{P1}{P2}, we find out that
\Bsist
  && \intQt |\dt\rhos|^2
  + \frac 12 \, \norma{\rhos(t)}_H^2 + \frac {\sigma}2 \, 
  \norma{\nabla \rhos(t)}_H^2 + \iO f_1(\rhos (t))
  \non
  \\
  && \leq c + c \iO |\rhos(t)|^2 + \frac 14 \intQt |\dt\rhos|^2
  + c \norma{\mus}_{\L\infty H}^2 . 
  \non
\Esist
On the other hand, by the chain rule and the Young inequality~\eqref{young} we have that
\Beq
  c \iO |\rhos(t)|^2 
  \leq c \iO |\rhozsig|^2 + \frac 14  \intQt |\dt\rhos|^2  
  + c\int_0^t \norma{\rhos(s)}_H^2\, ds.
  \non
\Eeq
Then, as $f_1$ is nonnegative, by accounting for \eqref{primastima}, 
with the help of the Gronwall lemma we infer that
\Beq
  \intQt |\dt\rhos|^2 + \norma{\rhos(t)}_H^2 + \sigma \, 
  \norma{\nabla \rhos(t)}_H^2 \leq c \quad \hbox{for all }\, t\in[0,T]. 
\non
\Eeq
Thus, we conclude~that
\Beq
  \norma{\rhos}_{\H1H}  + \sigma^{1/2} \norma{\rhos}_{\L\infty V} \leq c.
  \label{secondastima}
\Eeq

\step 
Third a priori estimate

We proceed formally and test \eqref{seconda} by $-\Delta\rhos$. Hence, 
integrating by parts and with respect to time, we deduce that 
\Bsist
  \hskip-1.5cm && \frac 12 \, \norma{\nabla \rhos(t)}_H^2
     + \sigma \intQt |\Delta \rhos|^2
    + \intQt \beta'(\rhos) |\nabla \rhos|^2
  \non
  \\
  \hskip-1.5cm && {}\leq \frac 12 \iO |\nabla\rhozsig|^2
   - \intQt \pi'(\rhos) |\nabla \rhos|^2
  + \intQt g'(\rhos) \nabla \mus \cdot \nabla \rhos 
  + \intQt g''(\rhos) \mus | \nabla \rhos |^2 ,
  \label{per3stima}
\Esist
where the equality  $\xis=\beta(\rhos)$ has been used along with the smoothness of $\beta$, 
according to our formal procedure. 
In fact, what is important is that the related term on the \lhs\ is nonnegative,~i.e., 
$$ 
  \intQt \beta'(\rhos) |\nabla \rhos|^2 \geq 0. 
$$
Concerning the \rhs\ of \eqref{per3stima}, we have that 
\ $\displaystyle \frac 12 \iO |\nabla\rhozsig|^2 \leq c $ \ 
due to \eqref{P1}, and the estimate 
$$
  {}- \intQt \pi'(\rhos) |\nabla \rhos|^2
  + \intQt g'(\rhos) \nabla \mus \cdot \nabla \rhos \leq
  c \int_0^t \norma{\nabla \rhos(s)}_H^2 ds + c \, \norma{\mus}_{\L2V}^2 
$$
owing to the boundedness of $\pi'$ and $g'$ (see
\accorpa{hpfdueg}{defbp}). About the last term,  \eqref{hpg} and \eqref{regmu} imply
$$
  \intQt g''(\rhos) \mus | \nabla \rhos |^2 \leq 0,
$$
so that the sign properties of $g''$ and $ \mus $ become crucial to control this term. 
Then, in view of~\eqref{primastima}, 
from \eqref{per3stima} it follows that 
$$
  \frac 12 \, \norma{\nabla \rhos(t)}_H^2
  + \sigma \intQt |\Delta \rhos|^2
  \leq c + c \int_0^t \norma{\nabla \rhos(s)}_H^2 \, ds
  \quad \hbox{for all }\, t\in[0,T], 
$$
and the Gronwall lemma and \eqref{secondastima} allow us to deduce that 
\Beq
  \norma{\rhos}_{\L\infty V}  + \sigma^{1/2} \norma{\rhos}_{\L2W} \leq c.
  \label{terzastima}
\Eeq
Note that here we have used the regularity theory for elliptic equations, 
owing to the bound on $\sigma \norma{\Delta \rhos}^2_2 $ 
and to the homogeneous Neumann boundary condition satisfied by~$\rhos $ (cf.~\eqref{regrho}). 
Finally, an easy consequence of  \eqref{secondastima}
and \eqref{terzastima} comes out from a comparison of terms in \eqref{seconda}, which yields 
\Beq
  \norma{\xis}_{\L2 H}  \leq c.
  \label{terzastimabis}
\Eeq

\step 
Fourth a priori estimate

As $\us = (  1+ 2 g(\rhos)) \mus $, by \eqref{hpfdueg} we have that 
\Bsist
  &|\us | \leq c \left( 1 + | \rhos |  \right) 
  | \mus |,
  & \non \\
  &\left| \nabla \us \right|  = \left| 2 
  g' (\rhos) \mus \nabla \rhos
  + (  1+ 2 g(\rhos)) \nabla \mus \right|
  {}\leq c \, | \mus |\, |\nabla \rhos |
  + c \left( 1 + | \rhos |  \right) |\nabla \mus | .
  &\non
\Esist
Now, taking \eqref{primastima} into account,
we see that $|\nabla\mus | $ is bounded in $\L2{\Lx2}$, 
while $|\mus| $~is bounded in $\L2{\Lx6}$ thanks to the Sobolev inequality~\eqref{sobolev}. 
On the other hand, \eqref{terzastima} 
provides a bound for $|\nabla \rhos| $ in $\L\infty{\Lx2} $ 
and for $| \rhos| $ in $\L\infty{\Lx6} $.  
Hence, using \holder\ inequality,
it is not difficult to check that the products  
$| \mus |\, |\nabla \rhos | $ and $ | \rhos | \, |\nabla \mus |$ are bounded in $\L2{\Lx{3/2}}$,
whereas $ | \rhos | \, | \mus |$ is even bounded in $\L2{\Lx{3}}$.
Therefore, we conclude that
\Beq
  \norma{\us}_{\L2{\Wx{1,3/2}}} \leq c \,.
  \label{quartastima}
\Eeq

\step
Fifth a priori estimate

Let us recall that \eqref{primastima} and \eqref{sobolev} imply 
the boundedness of $\{ \mus \} $ 
in the space   $\L\infty{\Lx{2}} \cap \L2{\Lx{6}}$. 
Then, we can apply \eqref{interpolx} with $p=2$, $q=6$, $\theta=1/2$, $r=3$ to see that
$$
  \norma{\mus (t) }_{3}^2 \leq \norma{\mus (t) }_{2} \norma{\mus (t) }_{6}
 \quad \hbox{for a.a. } \, t\in (0,T),
$$
whence squaring and integrating with respect to $t$ lead~to 
\Beq
  \label{boundL4L3} 
  \norma{\mus }^4_{{\L4{\Lx{3}}}} \leq  
  \norma{\mus }^2_{{\L\infty{\Lx{2}}}} \, 
  \norma{\mus }^2_{{\L2{\Lx{6}}}} \leq c.
\Eeq
Consider now \eqref{pri-variaz} which turns out to be 
a variational formulation 
of \eqref{prima}. As we want to prove that 
\Beq
  \norma{\dt \us}_{L^{4/3}(0,T;\Vp)} \leq c \,,
  \label{quintastima}
\Eeq
we use \eqref{pri-variaz} and let $v$ vary in $ L^{4}(0,T;V)$. 
By integrating with respect to time and invoking \eqref{hpkbis}, 
the boundedness of $g'$ and \holder's inequality, 
we obtain 
\Bsist
  &&\left|  \ioT  \< \dt \us (t), v (t) > \, dt  \right|
  \non
  \\
  && \leq \kappa^*  \norma{\nabla\mus}_{\L2H} \norma{\nabla v}_{\L2H} 
  +  c \int_0^T \norma{\mus (t)}_{3}  \norma{\dt\rhos (t)}_{2}
  \norma{v (t)}_{6} \, dt .
  \non
\Esist
Hence, in view of \eqref{primastima}, 
by applying the \holder\ and Sobolev inequalities (see~\eqref{sobolev}) in the time integral, 
we infer that 
\Bsist
  && \left|  \ioT  \< \dt \us (t), v (t) > \, dt  \right|
  \non
  \\
  && \leq c \norma{v}_{\L2V} +  c \norma{\mus}_{{\L4{\Lx{3}}}} 
  \norma{\dt\rhos}_{\L2H} \norma{v}_{L^{4}(0,T;V)} .
  \non
\Esist
Now, the continuous embedding $\L4V\subset\L2V$, 
\eqref{boundL4L3} and \eqref{secondastima} allow us to conclude that 
\Beq
  \left|  \ioT  \< \dt \us (t), v (t) > \, dt  \right| \leq c  \norma{v}_{L^{4}(0,T;V)} ,
  \non
\Eeq
whence \eqref{quintastima} follows.

\step
Passage to the limit

By the above estimates, there are a quadruplet $(\mu,\rho,\xi,u)$,
with $\mu\geq0$ \aeQ, and a function $k$ such~that
\accorpa{convmu}{convu} are satisfied as long~as 
\Beq
  \kamurhosig \to k  
  \quad \hbox{weakly star in $\LQ\infty$}
  \label{convkamurho}
\Eeq
at least for a subsequence $\tau=\tau_i{\scriptstyle\searrow}0$.
By the weak convergence of time derivatives,
the Cauchy conditions \eqref{cauchy} hold for the 
limit pair $(\rho,u)$. 
By \eqref{convrho}, \eqref{convu},
and the compact embedding~\eqref{compsobolev},
we can apply well-known strong compactness results
(see, e.g., \cite[Sect.~8, Cor.~4]{Simon})
and infer~that (possibly taking another subsequence)
\Bsist
  & \rhos \to \rho
  & \quad \hbox{strongly in $\C0{\Lx p}$ for $p<6$ and \aeQ}
  \label{strongrho}
  \\
  & \us \to u
  & \quad \hbox{strongly in $\L2{\Lx p}$ for $p<3$ and \aeQ}.
  \label{strongu}
\Esist
The weak convergence \eqref{convxi}, together with \eqref{strongrho} with $p=2$,
implies that $\xi\in\beta(\rho)$ \aeQ\ 
(see, e.g., \cite[Prop.~2.5, p.~27] {Brezis}), 
due to the maximal monotonicity of 
the operator induced by $\beta$ on $L^2(Q)$.
Now, we deal with the other nonlinear terms and the products.
We first observe that \eqref{strongrho}~also entails that
\Beq
  \phi(\rhos) \to \phi(\rho)
  \quad \hbox{strongly in $\C0{\Lx p}$ for $p<6$ and \aeQ} 
  \label{strongphirho}
\Eeq
for $\phi=g,g',\pi,1/(1+2g)$,
thanks to the \Lip\ continuity of such functions.
This is sufficient to establish equation~\eqref{secondez}.
Indeed, by accounting for~\eqref{convmu},
we see that the product $\mus g(\rhos)$ converges to $\mu g(\rho)$
weakly (e.g.) in~$\LQ2$.
On the other hand, \eqref{terzastima}~implies that
$\sigma\Delta\rhos$ converges to zero strongly in~$\LQ2$.
Now, we prove equations \accorpa{primau}{defu},
which involve the whole triplet $(\mu,\rho,u)$.
The first step is showing strong convergence for~$\mus$
and relation~\eqref{defu}.
By combining \eqref{strongphirho} with~\eqref{strongu},
we see~that
\Beq
  \mus = \frac \us {1+2g(\rhos)} \to \frac u {1+2g(\rho)}
  \quad \aeQ.
  \label{muaeQ}
\Eeq
This and \eqref{convmu} imply $\mu=u/(1+2g(\rho))$
and \eqref{defu} is proved.
Moreover, as $\graffe{\mus}$ is bounded in~$\LQ{10/3}$
by \eqref{primastima}, the Sobolev embedding $V\subset\Lx6$, and~\eqref{inter10-3-Q},
we can also deduce a strong convergence.
We \summariz e as follows:
\Beq
  \mus \to \mu
  \quad \hbox{strongly in $\LQ p$ for every $p<10/3$ \ and \aeQ}.
  \label{strongmu}
\Eeq
{}From this, we immediately infer that $\kamurhosig$ converges to $\kamurho$ \aeQ,
just by continuity.
Then, \eqref{convkamurho} implies $k=\kamurho$~and
\Beq
  \kamurhosig \to \kamurho
  \quad \hbox{strongly in $\LQ p$ for every $p<+\infty$}.
  \label{strongka}
\Eeq
Therefore, $\kamurhosig\nabla\mus$ converges to $\kamurho\nabla\mu$
weakly in $\LQ p$ for every $p<2$,
thanks to~\eqref{convmu},
and the choice $p=3/2$ yields
\Beq
  \intQ \kamurhosig\nabla\mus \cdot \nabla v  
  \to \intQ \kamurho\nabla\mu \cdot \nabla v
  \quad \hbox{for every $v\in\L3{\Wx{1,3}}$}.
  \non
\Eeq 
On the other hand, 
$\mus g'(\rhos)\dt\rhos$ converges to $\mu g'(\rho)\dt\rho$ 
weakly at least in~$\LQ1$, as one can easily see by combining
\eqref{convrho}, \eqref{strongphirho}, and~\eqref{strongmu}.
It follows that
\Beq
  \intQ \mus g'(\rhos) \dt\rhos \, v
  \to \intQ \mu g'(\rho) \dt\rho \, v
  \quad \hbox{for every $v\in\LQ\infty$}.
  \non
\Eeq
Moreover, \eqref{convu} holds.
Hence, we can conclude that
\Bsist
  && \ioT \< \dt u(t) , v(t) > \, dt
  + \intQ \kamurho\nabla\mu \cdot \nabla v
  = \intQ \mu g'(\rho) \dt\rho \, v
  \non
  \\
  && \quad \hbox{for every $v\in\L3{\Wx{1,3}}\cap\LQ\infty$}.
  \label{quasiprimau}
\Esist
Now, we observe that $\dt u\in\L{4/3}\Vp$ by \eqref{convu}
and that $\kamurho\nabla v\in\L2H$ by \eqref{convmu} and the boundedness of~$\kappa$. 
Finally, $\mu g'(\rho)\dt\rho\in\L{4/3}{\Lx{6/5}}$,
since $g'$ is bounded, $\dt\rho\in\L2H$,
and $\mu\in\L4{\Lx3}$ as a consequence of~\eqref{convmu},
$V\subset\Lx6$, and~\eqref{boundL4L3}).
Therefore, we can improve \eqref{quasiprimau} by a density argument
and see that the variational equation still holds for any $v\in\L4V$.
What we obtain is equivalent to~\eqref{primau},
and the proof is complete.


\section{Properties of the limit problem}
\label{Further}
\setcounter{equation}{0}

In this section, we prove Theorem~\ref{Linftyz}.
In the whole section, it is understood that the assumptions
of Theorem~\ref{Linftyz} are satisfied,
and sometimes we do not remind the reader about that.
As far as {\Gianni the first part of} Theorem~\ref{Linftyz} is concerned,
the true result regards ordinary variational inequalities
and we present it in the form of a lemma.
For convenience, we use the same notation $\rho$, etc.,
even though it is clear that everything is independent of~$x$:
{\Gianni the dot over the variable $\rho$ denotes the (time) derivative, here}.


\Blem \label{lem1}
Let \Perlimitatezza\ hold and $\rhomin\leq\rhoz\leq\rhomax$. Then
for every nonnegative function $\mu\in L^1(0,T)$, the differential inclusion
\Beq
  \dot\rho(t) + \beta(\rho(t)) + \pi(\rho(t)) - \mu(t) g'(\rho(t))
  \ni 0 
  \quad \aat
  \aand
  \rho(0) = \rhoz
  \label{eq:PJ8}  
\Eeq
has a unique solution $\rho\in W^{1,1}(0,T)$ such that
{\Gianni
\Beq
  \rhomin \leq \rho(t) \leq \rhomax  
  \aand
  \ximin \leq \xi(t) \leq \ximax \quad  \aat ,
  \label{eq:PJ9}
\Eeq
where
\Beq
  \xi(t)
  := - \bigl( \dot \rho(t) + \pi(\rho(t)) - \mu(t)g'(\rho(t)) \bigr)
  \in \beta(\rho(t)) .
  \non
\Eeq
}%
Moreover, there exists a constant $C>0$ such that if
$\mu_1, \mu_2\in L^1(0,T)$ and $\rhoz^1, \rhoz^2$ are two inputs and
$\rho_1, \rho_2$ are the corresponding solutions of \eqref{eq:PJ8}, then
for every $t \in [0,T]$ we have
\Bsist
  && {\Gianni |\rho_1-\rho_2|(t)}
  + \int_0^t |\dot\rho_1 - \dot\rho_2|(\tau)\,d\tau 
  \non
  \\
  && \leq C\left(|\rhoz^1- \rhoz^2| + \iot \bigl( {\GP(1+\mu_1)} |\rho_1 - \rho_2| + |\mu_1 - \mu_2| \bigr)(\tau)\,d\tau\right).
\label{eq:PJ9a}
\Esist

\Elem

\Bdim
The existence of a unique solution can easily be proved, 
{\GP e.g., by the iterated} Banach Contraction Principle,
{\Gianni due to the monotonicity of $\beta$ and to the \Lip\ continuity of the other nonlinearities}. 
In \eqref{eq:PJ9}, we only prove the upper inequalities since the proof
of the lower ones is quite similar. It suffices to prove the desired inequalities
for {\GP the solution $(\rho,\xi)$ of} the cut-off problem
{\GP
\Bsist
  && \dot\rho(t) + \xi(t) + \pst(\rho(t)) - \mu(t)\gst(\rho(t)) = 0 , \quad
  \xi(t) \in \beta(\rho(t))
  \quad \aat ,
  \qquad
  \label{eq:PJ10}
  \\
  && \rho(0) = \rhoz \,,
  \label{cauchyPG}
\Esist
}%
where $\pst$ and $\gst$ are defined~by
\Beq
  \pst(r) := \pi(\min\{r,\rhomax\})
  \aand
  \gst(r) := g'(\min\{r,\rhomax\}) \,.
  \non
\Eeq
We test \eqref{eq:PJ8} by $(\rho-\rhomax)^+$ and integrate.
{\Gianni
Recalling \accorpa{P13}{P15} and noting that
$\xi\geq\ximax$ and $g^*(\rho)=g'(\rhomax)$ where $\rho>\rhomax$,
we obtain}
\Bsist
   \frac 1 2 |(\rho(t)-\rhomax)^+|^2
  &\le& -\int_0^t \bigl( \xi - \ximax \bigr) (\rho-\rhomax)^+
  - \int_0^t \bigl(\ximax + \pi^*(\rhomax) \bigl) (\rho-\rhomax)^+
  \non
  \\
  && + \int_0^t \bigl( \pi(\rhomax) - \pi(\rho) \bigr) (\rho-\rhomax)^+
  + \int_0^t \mu \, g^*(\rho) (\rho-\rhomax)^+ \non
  \\
  &\le& \int_0^t \bigl( \pi(\rhomax) - \pi(\rho) \bigr) (\rho-\rhomax)^+
  \leq c \iot |(\rho-\rhomax)^+|^2
  \non
\Esist
and the assertion is obtained by the Gronwall argument. 
The second inequality follows
from the monotonicity of~$\beta$.
{\Gianni Moreover}, the lower bounds can be checked in a similar way.
To prove \eqref{eq:PJ9a}, we {\GP set}
$w_i(t) = \mu_i(t) g'(\rho_i(t)) - \pi(\rho_i(t))$,
$\xi_i(t) = w_i(t) -\dot\rho_i(t) $, $i=1,2$.
We have $(\xi_1-\xi_2)(\rho_1 - \rho_2) \ge 0$ {\Gianni almost everywhere}.
\vio{The function $\sign(\xi_1-\xi_2)$ 
{\GP(with $\sign(0)=0$)}
is~bounded and measurable, and so is 
$\sign(\rho_1-\rho_2)$. We now claim that by
testing the identity
\Beq
  (\xi_1-\xi_2) + (\dot\rho_1 - \dot\rho_2) = w_1 - w_2
  \label{prelip1}
\Eeq
by $\sign(\xi_1-\xi_2)$, we infer~that
\Beq
  |\xi_1-\xi_2| + \frac{d}{dt} |\rho_1 - \rho_2| \le |w_1 - w_2|
  \quad \hbox{a.e.\ in $(0,T)$}.
  \label{lip1}
\Eeq
Indeed, this is obvious for all $t$ such that $\sign(\xi_1-\xi_2)(t) = \sign(\rho_1-\rho_2)(t)$
or such that $\xi_1(t) = \xi_2(t)$. 
The remaining case is $\sign(\xi_1-\xi_2)(t) \ne 0$,
$\sign(\rho_1-\rho_2)(t) = 0$. 
For almost all $t$ with this property, we have
$\dot\rho_1(t) = \dot\rho_2(t)$, $\frac{d}{dt} |\rho_1 - \rho_2|(t) = 0$,
and \eqref{lip1} follows}.
{\Gianni
Using the Lipschitz continuity {\GP properties in \eqref{hpfdueg}}
and integrating \eqref{lip1} over~$(0,t)$, we obtain for~$t\in(0,T)$
\Beq
  \iot |\xi_1-\xi_2|(s) \, ds 
  + |\rho_1-\rho_2|(t)
  \leq c \left(|\rhoz^1- \rhoz^2| + \iot \bigl( {\GP(1+\mu_1)} |\rho_1 - \rho_2| + |\mu_1 - \mu_2| \bigr)(\tau)\,d\tau\right).
  \non
\Eeq
On the other hand, \eqref{prelip1} yields
\Beq
  \iot |\dot\rho_1 - \dot\rho_2|(s) \, ds
  \leq \iot \bigl( |w_1-w_2| + |\xi_1 - \xi_2| \bigr)(s) \, ds
  \non
\Eeq
and \eqref{eq:PJ9a} follows from the sum of the last two inequalities}.
\Edim

{\GP
Next, if $(\mu,\rho,\xi,u)$ is a solution to problem \Pblz, 
it is clear that, for almost all $x\in\Omega$,
the functions $\mu(x,\cdot)$ and $\rho(x,\cdot)$, and the constant $\rhoz(x)$
satisfy the assumptions of Lemma \ref{lem1}.
Thus, the first part of Theorem~\ref{Linftyz} concerning bounds \eqref{linftyz} is proved.
We derive an interesting consequence}.

\Bcor
\label{Boundedness}
Under the assumptions of Theorem~\ref{Linftyz},
let $(\mu,\rho,\xi,u)$ be a solution to problem \Pblz\
satisfying the regularity conditions specified in Theorem~\ref{Convergenza}.
Then 
\Beq
  \mu \in \LQ\infty
  \aand
  \dt\rho \in \LQ\infty .
  \label{boundedness}
\Eeq
\Ecor

\Bdim 
{\GP We already know that both $\xi$ and $\pi(\rho)$ are bounded.
Moreover}, $\mu g'(\rho)$
belongs to $\L\infty H\cap\L2{\Lx6}$ since $\mu$ does so and $g'(\rho)$ is bounded.
We see that also $\dt\rho$ belongs to such a space, just by comparison in~\eqref{secondez}.
It follows that $\dt\rho\in\L{7/3}{\Lx{14/3}}$ by~\eqref{inter7-14-3}. 
From this {\GP and assumption~\eqref{hpmuzbdd}}, we derive the boundedness of~$\mu$.
Indeed, we can reproduce the proof carried out in~\cite[Fifth a priori estimate]{CGPS7},
since that proof acts only on the equation for $\mu$
and works provided that an estimate of $\dt\rho$ in~$\L{7/3}{\Lx{14/3}}$ is known.
At this point, by comparing in~\eqref{secondez} once more,
we conclude that $\dt\rho$ is bounded as well.
\Edim

\Brem
\label{Linftys} 
The analogous estimate 
\Beq 
 \rhomin \leq \rhos \leq  \rhomax \quad \aeQ 
\label{linftys}
\Eeq
for the solution to problem \Pbl\ also holds provided that
\Beq
  \rhomin \leq \rhozsig \leq \rhomax
  \quad \aeO .
  \label{stimarhozs}
\Eeq
We prove one of the inequalities~\eqref{linftys}, the other one being similar.
We proceed as in the proof of Lemma \ref{lem1}, testing {\gianni\eqref{seconda}} by $(\rhos-\rhomax)^+$ and integrating.
By accounting for the second inequality \eqref{stimarhozs},
we easily obtain 
\Bsist
  && \frac 1 2 \iO |\rhosmrmp(t)|^2 
  + \sigma \intQt |\nabla(\rhos-\rhomax)^+|^2
  \non
  \\
  && \quad {}
  + \intQt \bigl( \xis - \ximax \bigl) \rhosmrmp  
  + \intQt \bigl(\ximax + \pi(\rhomax) \bigl) \rhosmrmp 
  \non
  \\
  && \leq \intQt \bigl( \pi(\rhomax) - \pi(\rhos) \bigl) \rhosmrmp
  + \intQt \mus \, g'(\rhos)  \rhosmrmp .
  \non
\Esist
Now, we observe that all the terms on the \lhs\ are nonnegative,
the third one thanks to \eqref{P13} and the monotonicity of~$\beta$
(as~before, the integrand vanishes whenever $\rhos\leq\rhomax$), 
the last one due to~\eqref{P14}. 
Concerning the \rhs, 
we show that the last integrand is nonpositive. 
Indeed, $g'$~is decreasing (see \eqref{hpg}), 
whence $g'(\rhos)\leq g'(\rhomax)\leq 0$ if $\rhos>\rhomax$,
and $\mus\geq0$.
By taking all this into account
and owing to the Lipschitz continuity of~$\pi$ (cf.~\eqref{defbp}),
we can apply the Gronwall lemma
and conclude that $\rhosmrmp=0$, i.e., $\rho\leq\rhomax$ \aeQ.
\Erem

\Brem
A sufficient condition for {\gianni\eqref{stimarhozs}} to hold at least for small~$\sigma$
is that $\rhozsig$ is given by~\eqref{P3}
and the hypotheses of Theorem~\ref{Linftyz} are reinforced by also assuming~that
\Beq
  \hbox{either} \qquad
  \infess\rhoz > \rhomin
  \enskip \hbox{and} \enskip
  \supess\rhoz < \rhomax 
  \qquad \hbox{or} \qquad
  \ximin \leq 0 \leq \ximax .
  \label{P12}
\Eeq
The proof is rather simple and we show just one of the desired inequalities
since the other one is quite similar.
We test {\gianni\eqref{P3}} by~$\rhozsmrmp$.
We easily obtain
\Bsist
  && \iO |\rhozsmrmp|^2
  + \sigma \iO |\nabla\rhozsmrmp|^2
  + \sigma \iO (\xizsig-\ximax) \rhozsmrmp
  \non
  \\
  && = \iO (\rhoz-\rhomax-\sigma\ximax) \rhozsmrmp . 
  \label{perstimarhozs}
\Esist
In the first case~\eqref{P12}, we set $\delta:=\rhomax-\supess\rhoz$
and take $\sigma^*>0$ such that $\sigma^*\,|\ximax|\leq\delta$.
Then, for $\sigma\leq\sigma^*$, we have 
$\rhoz-\rhomax-\sigma\ximax \leq -\delta+\sigma^*|\ximax| \leq 0$ \aeO,
so that the \rhs\ of \eqref{perstimarhozs} is nonpositive.
In the second case~\eqref{P12}, the same conclusion trivially holds.
As the last two terms on the \lhs\ are nonnegative
(since \eqref{P13} holds, $\beta$~is monotone, 
and the third integrand vanishes whenever $\rhozsig\leq\rhomax$),
we~conclude that $\rhozsmrmp=0$, whence $\rhozsig\leq\rhomax$.
\Erem

\step 
Proof of the {\Gianni second part of} Theorem~\ref{Linftyz}

{\Gianni Assume thus} that $\kappa(\mu, \rho) = \kappa_0$
and set for simplicity $\kappa_0 = 1$. 
The system now reads
\Bsist
  & \< \dt u(t), v > + \displaystyle\iO \nabla\mu(t) \cdot \nabla v
  = \displaystyle\iO \mu \, g'(\rho) \, \dt\rho \, v \qquad
  &
  \non
  \\
  &\hskip6cm \hbox{for all } v\in V \hbox{ and a.a. } t\in (0,T),
  &
  \label{pr1}
  \\
   & u= (1+2g(\rho))\mu 
  \quad \aeQ,&\label{de1}
  \\
   &   \dt\rho + \xi + \pi(\rho)
   = \mu \, g'(\rho)
  \aand \xi \in \beta(\rho)
  \quad \aeQ, & \label{se1}
  \\
  & \mu(0) = \muz
  \aand
  \rho(0) = \rhoz
   \quad \aeO.& 
  \label{ca1}
\Esist
Let {\Gianni $(\mu_i,\rho_i,\xi_i,u_i)$}, $i=1,2$ be two solutions of \eqref{pr1}--\eqref{ca1}.
We integrate \eqref{pr1} in time from $0$ to $t$ and subtract the equation
with index $2$ from the one with index $1$. 
We test the result by {\Gianni $v = (\mu_1 - \mu_2)(t)$}
and obtain, by virtue of Corollary \ref{Boundedness}, that
\Bsist
  \non
  && \int_\Omega (u_1 - u_2)(\mu_1 - \mu_2)(t)
  + {\Gianni \frac12 \frac{d}{dt} \iO \left|\iot  \nabla(\mu_1 - \mu_2)d\tau \right|^2 }
  \\
  \label{uni1}
  &&\le {\Gianni c} \int_\Omega \left(|\mu_1 - \mu_2|(t) \int_0^t \left(|\mu_1 - \mu_2|
  + |\rho_1 - \rho_2| + |\dt\rho_1 - \dt\rho_2|\right)(\tau) d\tau\right).
\Esist
{\GP
In addition, from Lemma~\ref{lem1} (see, in particular, \eqref{eq:PJ9a}) 
and \holder's inequality it follows that
\Bsist
  && \iO \left(\int_0^t |\dt\rho_1 - \dt\rho_2|(\tau)\,d\tau\right)^2 \le
  c \int_\Omega\left(\int_0^t(|\rho_1 - \rho_2|
  + |\mu_1 - \mu_2|)(\tau)\,d\tau\right)^2 ,
  \qquad
  \label{uni2}
  \\
  && \iO |\rho_1 - \rho_2|^2(s) \le
  D \int_0^s \!\!\iO \bigl( |\rho_1 - \rho_2|^2 + |\mu_1 - \mu_2|^2 \bigr)(\tau)\,d\tau
\label{uni2a}
\Esist
for every $t,s\in[0,T]$,
thanks to the boundedness for $\mu_1$ ensured by Corollary~\ref{Boundedness}}.
{\GP
Note that the constant $D$ in \eqref{uni2a} is marked for later reference.}

{\GP Now, we observe that the inequalities
\Beq
  (u_1 - u_2) (\mu_1 - \mu_2)
  \geq |\mu_1 - \mu_2|^2 
  - 2 \mu_1 \bigl( g(\rho_1) - g(\rho_2) \bigr) (\mu_1 - \mu_2)
  \geq \frac 12 \, |\mu_1 - \mu_2|^2
  - c |\rho_1 - \rho_2|^2
  \non
\Eeq
hold \aeQ.
Thus, by integrating}
\vio{\eqref{uni1} from $0$ to~$s$, $s \in (0,T)$, 
{\GP and ignoring a positive term on the \lhs, we obtain}
\Bsist
  \non
  &&\int_0^s\!\!\int_\Omega |\mu_1 - \mu_2|^2(t)\, dt
  \le c\int_0^s\!\!\int_\Omega |\rho_1 - \rho_2|^2(t)\, dt
  + c\left(\int_0^s\!\!\int_\Omega |\mu_1 - \mu_2|^2(t)\, dt\right)^{1/2}
  \\
  \label{uni1a}
  && \times
  \left(\int_0^s\!\!\int_\Omega \left(\int_0^t \left(|\mu_1 - \mu_2|
  + |\rho_1 - \rho_2| + |\dt\rho_1 - \dt\rho_2|\right)(\tau) d\tau\right)^2 dt\right)^{1/2}.
\Esist
{\GP Hence}, using Young's inequality and \eqref{uni2}, {\GP we have~that}
\Bsist \non
  &&\int_0^s\!\!\int_\Omega |\mu_1 - \mu_2|^2(t)\, dt
  \le c\int_0^s\!\!\int_\Omega |\rho_1 - \rho_2|^2(t)\, dt\\
  \label{uni1b}
  &&
  + c\int_0^s\!\!\int_\Omega \left(\int_0^t \left(|\mu_1 - \mu_2|
  + |\rho_1 - \rho_2|\right)(\tau) d\tau\right)^2 dt.
\Esist
We now {\GP multiply \eqref{uni1b} by $2D$ and add it to~\eqref{uni2a}.
Thus, we} obtain an inequality of the form
$\Phi(s) \le c \int_0^s \Phi(t) dt$, with
$$
\Phi(s) = \int_\Omega |\rho_1 - \rho_2|^2(s) + \int_0^s\!\!\int_\Omega |\mu_1 - \mu_2|^2(t)\, dt.
$$
From the Gronwall argument, {\GP it is \sfw\ to deduce that} $\Phi(s) = 0$ for all $s$, hence,
$\mu_1 = \mu_2$, $\rho_1 = \rho_2$,
which implies uniqueness.}

The $L^2$ bound for $\dt\mu$ can be established in the following way.
Assume first that $\muz\in W$. 
We extend $\mu$ by $\mu_0$ and $\rho$ by $\rho_0$ for $t<0$. 
{\GP Then, equation} \eqref{pr1} then can be written~as
\Beq
  \< \dt u(t), v > + \iO \nabla\mu(t) \cdot \nabla v
  = \iO \psi(t) \, v 
  {\GP \quad \hbox{for all $v\in V$ and a.a.\ $t\in(0,T)$,}}
  \label{reg1}
\Eeq
{\GP where $\psi$ is defined by
$\psi(t) = \bigl(\mu g'(\rho)\dt\rho\bigr)(t)$ for $t>0$
and $\psi(t)=-\Delta\muz$ for $t<0$.
We observe that $\psi\in L^\infty(-T,T;H)$ thanks to Corollary~\ref{Boundedness}
and to our assumption on~$\muz$.
Next, we integrate \eqref{reg1} in time from $(t{-}h)$ to $t$ for any fixed $t \in (0,T)$
and a small~{\GP$h>0$}, with the intention to let $h$ tend to zero, 
and test the resulting equality by $\mu(t)-\mu(t-h)$.
We obtain
\Bsist
  && \iO \bigl( u(t) - u(t-h) \bigr) \bigl( \mu(t) - \mu(t-h) \bigr) +
  {\GP \frac12 \iO \frac{d}{dt}
    \left| \int_{t-h}^t \nabla\mu(\tau) \, d\tau \right|^2}
  \non
  \\
  && {\GP = \iO \left( \int_{t-h}^t \psi(\tau) \, d\tau \right) \, \bigl( \mu(t) - \mu(t-h) \bigr)} 
  \non
  \\
  && {\GP \leq \frac 14 \iO |\mu(t) - \mu(t-h)|^2
  + \left\| \int_{t-h}^t \psi(\tau) \, d\tau \right\|_H^2 }
  \non
  \\
  && {\GP \leq \frac 14 \iO |\mu(t) - \mu(t-h)|^2
  + c \, h^2 }
  \label{reg2}
\Esist
{\GP Now, we recall that \eqref{de1} holds, that $g$ is nonnegative and \Lip\ continuous,
and that $\mu$ and $\dt\rho$ are bounded by Corollary~\ref{Boundedness}. 
Hence, we easily derive~that
\Bsist
  && \bigl( u(t) - u(t-h) \bigr) \bigl( \mu(t) - \mu(t-h) \bigr)
  \non
  \\
  && \geq |\mu(t) - \mu(t-h)|^2
  - 2 \mu(t) \, |g(\rho(t)) - g(\rho(t-h))| \, |\mu(t) - \mu(t-h)| 
  \non
  \\
  && \geq |\mu(t) - \mu(t-h)|^2
  - c \, h \, |\mu(t) - \mu(t-h)| 
  \geq \frac 12 \, |\mu(t) - \mu(t-h)|^2
  - c \, h^2 .
  \non
\Esist
Therefore, by integrating \eqref{reg2} from $0$ to~$T$,
forgetting the nonnegative term that involves~$\nabla\mu$, and rearranging,
we obtain
\Beq
  \ioT \!\! \iO \left| \mu(t) - \mu(t-h) \right|^2 dt
  \leq c \, h^2 + c \iO \left| \int_{-h}^0 \nabla\muz \, d\tau \right|^2
  \leq c \, h^2 .
  \non
\Eeq
As $h>0$ is arbitrarily small, this implies that $\dt\mu\in\LQ2$.
At this point, we are allowed to use the Leibniz rule
for the time derivative~$\dt u$;
then, from \accorpa{pr1}{de1} we infer that
the equation
\Beq
  \bigl( \coeff \bigr) \dt\mu + \mu g'(\rho) \dt\rho 
  - \Delta\mu
  = 0 
  \label{strongform}
\Eeq
holds at least in the sense of distributions}.
By comparison,  {\GP we deduce that $\Delta \mu \in L^2(Q)$, whence $\mu\in\L2W$}. 
Using the identity
$$
  -\int_\Omega \dt\mu\,\Delta\mu
  = \frac12\frac{d}{dt} {\GP \iO |\nabla\mu|^2}
  \quad \hbox{a.e.\ in $(0,T)$},
$$
we see that $\nabla\mu \in L^\infty(0,T; L^2(\Omega))$.
{\GP Thus, the regularity \eqref{piuregolare} is established if $\muz\in W$.}

Let now {\GP $\mu_0 \in V\cap\Linfty$} be arbitrary, and 
{\GP consider a sequence $\{\mu^0_k\}\subset W$ bounded in $\Linfty$ and converging to $\muz$
in~$V$} as $k \to \infty$. 
Let {\GP $(\mu_k,\rho_k,\xi_k,u_k)$}
be the corresponding solutions to \eqref{pr1}--\eqref{ca1}.
{\GP Then, we can use equation~\eqref{strongform}
written with the index~$k$
and test it by $\dt\mu_k$.
We obtain}
\Beq
  \label{reg4}
  {\GP \iO |\dt\mu_k(t)|^2}
  + \frac12\frac{d}{dt} {\GP \iO |\nabla\mu_k(t)|^2}
  \le \int_\Omega |\psi_k(t)| \, |\dt\mu_k(t)| ,
\Eeq
with {\GP an obvious choice of $\psi_k\in\LQ2$} bounded in this space
{\GP(even better)}
independently of~$k$. 
{\GP By time integration, it is \sfw\ to}
obtain a bound for $\|\dt\mu_k\|_{L^2(Q)}$ {\GP and for $\norma{\nabla\mu_k}_{\L\infty H}$}
independent of~$k$. 
{\GP
Then, by weak star compactness we infer~that
\Beq
  \mu_k \to \tilde\mu
  \quad \hbox{weakly star in $\H1H\cap\L\infty V$}
  \non
\Eeq
at least for a subsequence, which implies (see, e.g., \cite[Cor.~4, p.~85]{Simon})
strong convergence in $\C0H$.
In particular, $\tilde\mu(0)=\muz$.
On the other hand, $(\mu_k,\rho_k,\xi_k,u_k)$ satisfies 
the estimates stated in Lemma~\ref{lem1}
and the boundedness properties for $\mu_k$ and $\dt\rho_k$ given by Corollary~\ref{Boundedness},
which are uniform with respect to~$k$.
This yields weak or weak star limits $\tilde\rho$ and~$\tilde\xi$.
Moreover, strong convergence in $L^1(Q)$ for $\{\rho_k\}$ and $\{\dt\rho_k\}$ is ensured 
via a Cauchy sequence argument based on~\eqref{eq:PJ9a}, integration over $\Omega$, and Gronwall's lemma. Hence, $\{\mu_k\}$, $\{\rho_k\}$, $\{\dt\rho_k\}$ converge strongly in  $L^p(Q)$ for every $p \in [1,\infty)$.  
At this point, it is not difficult to verify that $(\tilde\mu,\tilde\rho,\tilde\xi,\tilde u)$, with the corresponding~$\tilde u$, actually solves problem~\Pblz\ and thus coincides with the unique solution $(\mu,\rho,\xi,u)$.
Therefore, the proof is complete.}}

\section*{Acknowledgments}
{\GP The authors gratefully 
acknowledge the warm hospitality of the IMATI of CNR in Pavia, 
the Institute of Mathematics of the Czech Academy of Sciences 
in Prague, and the WIAS in Berlin. The present paper benefits 
from the GA\v CR Grant P201/10/2315 and RVO:~67985840 for~PK,
the MIUR-PRIN Grant 2010A2TFX2 ``Calculus of variations'' for PC and~GG, 
and the FP7-IDEAS-ERC-StG Grant \#200947 (BioSMA) for PC and~JS.
The work of JS was also supported by the DFG Research Center 
{\sc Matheon} in Berlin.}



\vspace{3truemm}

\Begin{thebibliography}{10}

\bibitem{Barbu}
V. Barbu,
``Nonlinear semigroups and differential equations in Banach spaces'',
Noord\-hoff,
Leyden,
1976.


\bibitem{Brezis}
H. Brezis,
``Op\'erateurs maximaux monotones et semi-groupes de contractions
dans les espaces de Hilbert'',
North-Holland Math. Stud.
{\bf 5},
North-Holland,
Amsterdam,
1973.


\bibitem{CGPS3} 
P. Colli, G. Gilardi, P. Podio-Guidugli, J. Sprekels,
Well-posedness and long-time behaviour for a nonstandard viscous Cahn-Hilliard 
system, {\it SIAM J. Appl. Math.} {\bf 71} (2011) 1849-1870.

\bibitem{CGPS4} 
P. Colli, G. Gilardi, P. Podio-Guidugli, J. Sprekels,
An asymptotic analysis for a nonstandard Cahn-Hilliard system with viscosity, 
{\GP{\it Discrete Contin. Dyn. Syst. Ser. S} {\bf 6} (2013) 353-368.}

\bibitem{CGPS6} 
P. Colli, G. Gilardi, P. Podio-Guidugli, J. Sprekels,
Global existence and uniqueness for a singular/degenerate  
Cahn-Hilliard system with viscosity, preprint 
WIAS-Berlin n.~1713 (2012), pp.~1-28.

\bibitem{CGPS7} 
P. Colli, G. Gilardi, P. Podio-Guidugli, J. Sprekels,
Global existence for a strongly coupled 
Cahn-Hilliard system with viscosity, 
{\it Boll. Unione Mat. Ital. (9)} {\bf 5} (2012) 495-513.

\bibitem{CS1}
P. Colli, J. Sprekels, On a Penrose--Fife model
with zero interfacial energy leading to a phase-field system of
relaxed Stefan type, {\it Ann. Mat. Pura Appl. (4)} {\bf 169} (1995)
269-289.

\bibitem{CS2} 
P. Colli, J. Sprekels, Global solution to the
Penrose-Fife phase-field model with zero interfacial energy and
Fourier law, {\it Adv. Math. Sci. Appl.} {\bf 9} (1999) 383-391.


\bibitem{GKS}
G. Gilardi, P. Krej{\v{c}}{\'{\i}}, J. Sprekels,
{Hysteresis in phase-field models with thermal memory},
{\it Math. Methods Appl. Sci.}
{\bf 23}
(2000)
{909-922}.
              
\bibitem{KS1}
P. Krej{\v{c}}{\'{\i}}, J. Sprekels,
{Hysteresis operators in phase-field models of {P}enrose-{F}ife type},
{\it Appl. Math.}
{\bf 43}
(1998)
{207-222}.

\bibitem{KS2}
P. Krej{\v{c}}{\'{\i}}, J. Sprekels,
{A hysteresis approach to phase-field models},
{\it Nonlinear Anal.}
{\bf 39}
(2000)
{569--586}.
       
\bibitem{KS3}
P. Krej{\v{c}}{\'{\i}}, J. Sprekels,
{Phase-field models with hysteresis},
{\it J. Math. Anal. Appl.}
{\bf 252}
(2000)
{198-219}.

\bibitem{KSZ}
P. Krej{\v{c}}{\'{\i}}, J. Sprekels, S. Zheng,
{Asymptotic behaviour for a phase-field system with hysteresis},
{\it J. Differential Equations}
{\bf 175} (2001)
{88-107}.





\bibitem{Podio}
P. Podio-Guidugli, 
Models of phase segregation and diffusion of atomic species on a lattice,
{\it Ric. Mat.} {\bf 55} (2006) 105-118.


\bibitem{Simon}
J. Simon,
{Compact sets in the space $L^p(0,T; B)$},
{\it Ann. Mat. Pura Appl.~(4)} {\bf 146} (1987) 65--96.

\End{thebibliography}

\End{document}

\bye